\documentclass[proceedings]{aofa}
\usepackage{amsmath}
\usepackage{graphicx}
\usepackage{url}
\usepackage{cite}

\usepackage[utf8]{inputenc}
\usepackage{subfigure}

\newtheorem{theo}{Theorem}

\newtheorem{lemma}[theo]{Lemma}

\newtheorem{rem}{Remark}

\newcommand{\E}{\mathbb{E}}

\newcommand{\Oh}{\mathcal{O}}

\allowdisplaybreaks

\title{Additive functionals of $d$-ary increasing trees}

\author[Dimbinaina Ralaivaosaona \and Stephan Wagner]{Dimbinaina Ralaivaosaona\addressmark{1}\and Stephan Wagner\addressmark{1}\thanks{This material is based upon work supported financially by the National Research Foundation of South Africa under grant number 96236.}}
\address{Department of Mathematical Sciences, Stellenbosch University, Private Bag X1, Matieland 7602, South
  Africa, \url{{naina,swagner}@sun.ac.za}}


\begin{document}
\maketitle

\begin{abstract}
A tree functional is called additive if it satisfies a recursion of the form $F(T) = \sum_{j=1}^k F(B_j) + f(T)$, where $B_1,\ldots,B_k$ are the branches of the tree $T$ and $f(T)$ is a toll function. We prove a general central limit theorem for additive functionals of $d$-ary increasing trees under suitable assumptions on the toll function. The same method also applies to generalised plane-oriented increasing trees (GPORTs).  One of our main applications is a log-normal law that we prove for the size of the automorphism group of $d$-ary increasing trees, but many other examples (old and new) are covered as well.
\end{abstract}

\section{Introduction}\label{sec:intro}

In this paper, we are interested in functionals of rooted trees that satisfy an \emph{additive} relation, i.e.~a recursion of the form
\begin{equation}\label{eq:F_recursion}
F(T) = \sum_{j=1}^k F(B_j) + f(T),
\end{equation}
where $B_1,\ldots,B_k$ are the branches of the tree $T$ and $f(T)$ is a so-called toll function, which often only depends on specific features of the tree such as the size or the root degree, but can in principle be arbitrary. The trees in our context will be labelled; it is assumed that the toll function only depends on the relative order of the labels, not the labels themselves, so that it is also well-defined if the labels are not necessarily $1,2,\ldots,n$. It is consistent with~\eqref{eq:F_recursion} to assume that we have the identity $F(\bigodot) = f(\bigodot)$ for the tree $T= \bigodot$ consisting only of a single labelled vertex. Important examples include
\begin{itemize}
\item the number of leaves, which corresponds to the toll function 
$$f(T) = \begin{cases} 1 & |T| = 1, \\ 0 & \text{otherwise}.\end{cases}$$
\item the number of vertices of outdegree $k$, in which case one can simply take
$$f(T) = \begin{cases} 1 & \text{if the root of $T$ has outdegree $k$,} \\ 0 & \text{otherwise}.\end{cases}$$
\item the internal path length, i.e., the sum of the distances from the root to all vertices, which can be obtained from the toll function $f(T) = |T|-1$.
\item the log-product of the subtree sizes \cite{meir1998logproduct}, also called the ``shape functional'' \cite{fill2004limiting}, corresponding to $f(T) = \log |T|$,
\item the logarithm of the size of the automorphism group: here, it is not difficult to see that the toll function is $f(T) = \log(R(T))$, where $R(T)$ is the size of the symmetry group of the collection of root branches.
\end{itemize}
Such functionals also arise frequently in the study of divide-and-conquer algorithms, e.g. quicksort \cite{hwang2002phase}. An alternative viewpoint is based on the notion of \emph{fringe subtrees}: a fringe subtree of a tree is a subtree induced by a vertex and all its descendants. If we let $\mathcal{F}(T)$ denote the collection of all fringe subtrees of a tree $T$, then it is easy to verify that
$$F(T) = \sum_{S \in \mathcal{F}(T)} f(S).$$
In particular, the number of occurrences of a specific tree as a fringe subtree is an additive functional (corresponding to the case that the toll function $f$ is an indicator function), and every additive functional can be obtained as a linear combination of such special functionals.

There are several recent papers providing central limit theorems for rather general additive tree functionals \cite{fill2004limiting,devroye2002limit,fill2005singularity,holmgren2015limit,janson2016asymptotic,wagner2015central}. Specifically, Holmgren and Janson \cite{holmgren2015limit} proved such a central limit theorem for binary increasing trees (which are also equivalent to binary search trees) and recursive trees. Both are instances of so-called \emph{increasing trees}: labelled trees with the additional property that the labels increase along any path starting at the root.

Varieties of increasing trees were studied systematically in \cite{bergeron1992varieties} (see also \cite[Section 1.3.3]{drmota2009random}). The exponential generating function $Y(x)$ associated with a variety of increasing trees satisfies a differential equation of the characteristic shape
\begin{equation}\label{eq:diffeq1}
Y'(x) = \Phi(Y(x)),\qquad Y(0) = 0
\end{equation}
for some function $\Phi(t)$. Varieties of increasing trees for which a uniformly random tree with a given number of vertices can also be generated by a growth process have been of particular interest. There are three such types \cite{panholzer2007level}:
\begin{itemize}
\item The variety of recursive trees is perhaps the most basic instance: these are simply labelled rooted unordered trees (``unordered'' meaning that the order of branches does not matter) with the aforementioned property that the labels increase along paths starting at the root. Uniformly random recursive trees can be obtained by the following growth process: starting from a single vertex (the root, carrying label $1$), the vertex labelled $n$ is attached in the $n$-th step to one of the previous vertices, chosen uniformly at random. As mentioned earlier, the order of children attached to a vertex does not matter. To obtain a canonical representation, one can e.g. always make the newly added vertex the rightmost child. In this example, the function $\Phi$ is the exponential function. It is easy to see that the generating function is $Y(x) = -\log(1-x)$, and for every positive integer $n$, there are $(n-1)!$ recursive trees.
\item Plane-oriented recursive trees (PORTs) differ from recursive trees in only one aspect: trees are regarded as embedded in the plane, the order of branches is taken into account. The growth process to generate uniformly random PORTs follows a ``preferential attachment'' rule: it is essentially the same as for recursive trees, but the probability that the vertex labelled $n$ is attached to a specific vertex $v$ is proportional to $1$ plus the current outdegree of $v$. Here, we have $\Phi(t) = (1-t)^{-1}$, so the generating function is $Y(x) = 1 - \sqrt{1-2x}$, and the number of plane oriented recursive trees with $n$ vertices is $(2n-3)!!$.

Generalised plane oriented recursive trees (GPORTs) are obtained by introducing an additional parameter: for some positive real number $\alpha$, we let the probability that the vertex labelled $n$ is attached to a specific vertex $v$ be proportional to $\alpha$ plus the current outdegree of $v$. An equivalent description uses weighted PORTs: to each PORT $T$, we associate a weight based on its outdegrees. If $N_j(T)$ is the number of vertices whose outdegree is $j$, we set
$$w(T) = \prod_{j \geq 1} \binom{\alpha+j-1}{j}^{N_j(T)}.$$
In choosing a random GPORT, the probability of a tree to be chosen is proportional to its weight. In the exponential generating function $Y(x)$, each tree is also weighted with $w(T)$. The function $\Phi$ in~\eqref{eq:diffeq1} is now given by $\Phi(t) = (1-t)^{-\alpha}$. It follows that $Y(x) = 1- (1-(\alpha+1)x)^{1/(1+\alpha)}$, the total weight of all trees with $n$ vertices is $\prod_{j=1}^{n-1} ((\alpha+1)j-1)$.
\item Finally, we have the variety of $d$-ary increasing trees, which will be the focus of this paper: here, every vertex has $d$ possible places to which a child can be attached (for example, in the binary case, there are left and right children). In the construction of uniform $d$-ary increasing trees by a growth process, we simply attach the vertex labelled $n$ to one of the $(d-1)(n-1)+1$ places available in total, once again selected uniformly at random. Therefore, the probability that the new vertex is attached to an existing vertex $v$ is proportional to $d$ minus the current outdegree of $v$ (in particular, if $v$ already has $d$ children, no further vertices can be attached to it). Here, $\Phi(t) = (1+t)^d$ and $Y(x) = (1-(d-1)x)^{-1/(d-1)} - 1$. The total number of $d$-ary increasing trees with $n$ vertices is
$\prod_{j=1}^{n-1} ((d-1)j+1)$.
\end{itemize}
 
\begin{rem}
We remark that recursive trees and $d$-ary increasing trees can also be seen as weighted PORTs, with weights
$$w(T) = \prod_{j \geq 1} c_j^{N_j(T)},$$
where $c_j = \frac{1}{j!}$ for recursive trees (to factor out the different ways of ordering the branches) and $c_j = \binom{d}{j}$ (to take the $d$ possible points of attachment into account) respectively.
\end{rem}

In the following, we state and prove a central limit theorem for additive tree functionals of uniformly random $d$-ary increasing trees under certain technical conditions on the toll function. As mentioned earlier, binary increasing trees (as well as recursive trees) have already been covered in \cite{holmgren2015limit}. Since the approach in \cite{holmgren2015limit} is based on representations of binary increasing trees and recursive trees that are not available for other classes of increasing trees, and the generating function method of \cite{wagner2015central} requires the resulting differential equations to be explicitly solvable, which is also not the case, we use a different approach based on moments, as in a paper of Fuchs \cite{fuchs2012limit} on the number of fringe subtrees of given size (which is also an additive functional). Although we only discuss the case of $d$-ary increasing trees in (some) detail, our method also applies to GPORTs, for which we only state the corresponding result in the following section.

This extended abstract only summarises the proof our main theorem and lists some interesting examples to which our result can be applied. Technical details and proofs of all intermediate lemmas will be provided in the full version of this paper.  

\section{The general central limit theorem}

Let us now formulate our main result. In the following, $d$ is fixed, and $T_n$ always denotes a random $d$-ary increasing tree of order $n$ (except for Theorem~\ref{thm:gport}). We assume that the toll function $f(T)$ satisfies the following conditions:
\begin{itemize}
\item[(C1)] $f(T)$ is bounded,
\\
\item[(C2)] $\displaystyle \sum_{k\geq 1}\frac{\mathbb{E}|f(T_k)|}{k}<\infty$ and $\mathbb{E}|f(T_n)| \to 0$ as $n\to\infty$. 
\end{itemize}

Under these assumptions, our central limit theorem for additive functionals reads as follows:

\begin{theo}\label{thm:main}
Let $T_n$ be a uniformly random $d$-ary increasing tree with $n$ vertices. If the toll function $f(T)$ satisfies (C1) and (C2), then there exist constants $\mu$ and $\sigma$ such that the mean and variance of $F(T_n)$ are asymptotically
$$\E(F(T_n)) = \mu n + \frac{\mu}{d-1}+o(1), \, \text{ \  } \,\mathrm{Var} (F(T_n))=\sigma^2n+o(n).$$ 
The constants $\mu$ and $\sigma$ can be represented as
\begin{equation}\label{eq:mu-eq}
\mu = (d-1) \sum_{T} f(T) \prod_{j=1}^{|T|} \frac{1}{(d-1)j+d}
\end{equation}
and
\begin{align*}
\sigma^2 =
& -\frac{\mu^2}{d-1}-(d-1)\sum_{T}\frac{f(T)^2-2f(T)(F(T)-\mu|T|)}{\prod_{j=1}^{|T|} ((d-1)j+d)} +\\
& d\sum_{T_1}\sum_{T_2}\frac{(d-1)^{1-|T_1|-|T_2|}f(T_1)f(T_2)}{(|T_1|-1)!(|T_2|-1)!}\int_0^1 \phi_{|T_1|}(x)\phi_{|T_2|}(x)dx,
\end{align*}
where 
$$
\phi_k(x)=(1-x)^{-1}\int_{x}^1(1-w)^{d/(d-1)}w^{k-1}dw.
$$
The sums are taken over all $d$-ary increasing trees.
If $\sigma \neq 0$, then the renormalised random variable $(F(T_n)-\mu n)/\sqrt{\sigma^2 n}$
converges weakly to a standard normal distribution.
\end{theo}

\begin{rem}\label{rem:constant}
We remark that the result remains true if conditions (C1) and (C2) hold for a shifted version $f(T) + c$ ($c$ any constant) of the toll function rather than the toll function itself, since this changes $F(T)$ only by the deterministic quantity $c|T|$.
\end{rem}

\begin{rem}
By means of the Cram\'er-Wold device, we also obtain joint normal distribution of tuples of additive functionals.
\end{rem}

As mentioned earlier, the method used in proving Theorem~\ref{thm:main} also applies to GPORTs. Without going into detail, let us just state the corresponding theorem:

\begin{theo}\label{thm:gport}
Let $T_n$ be a random GPORT (with fixed parameter $\alpha$) with $n$ vertices. If the toll function $f(T)$ satisfies (C1) and (C2), then there exist constants $\mu$ and $\sigma$ such that the mean and variance of $F(T_n)$ are asymptotically
$$\E(F(T_n)) = \mu n - \frac{\mu}{\alpha+1} +o(1), \, \text{ \  } \,\mathrm{Var} (F(T_n))=\sigma^2n+o(n).$$ 
The constants $\mu$ and $\sigma$ can be represented as
$$
\mu = (\alpha+1) \sum_{T} w(T) f(T) \prod_{j=1}^{|T|} \frac{1}{(\alpha+1)j+\alpha}
$$
and
\begin{align*}
\sigma^2 =
& \frac{\mu^2}{\alpha+1}-(\alpha+1)\sum_{T} w(T) \frac{f(T)^2-2f(T)(F(T)-\mu|T|)}{\prod_{j=1}^{|T|} ((\alpha+1)j+\alpha)} +\\
& \alpha \sum_{T_1}\sum_{T_2} w(T_1)w(T_2) \frac{(\alpha+1)^{1-|T_1|-|T_2|}f(T_1)f(T_2)}{(|T_1|-1)!(|T_2|-1)!}\int_0^1 \varphi_{|T_1|}(x)\varphi_{|T_2|}(x)dx,
\end{align*}
where
$$
\varphi_k(x)=\int_{x}^{1}(1-w)^{\alpha/(\alpha+1)}w^{k-1}dw.
$$
The sums are taken over all PORTs, weighted by $w(T)$.
If $\sigma \neq 0$, then the renormalised random variable $(F(T_n)-\mu n)/\sqrt{\sigma^2 n}$ converges weakly to a standard normal distribution.
\end{theo}

\section{Preliminaries}

Recall that the exponential generating function $Y(x)$ of $d$-ary increasing trees satisfies the differential equation
\begin{equation}\label{eq:diffeq2}
Y'(x) = \Phi(Y(x)),\qquad Y(0) = 0,
\end{equation}
where $\Phi(t) = (1+t)^d$. The explicit solution is given by $Y(x) = (1-(d-1)x)^{-1/(d-1)} - 1$, and the total number of $d$-ary increasing trees with $n$ vertices is
$$Y_n = n! \cdot [x^n] Y(x) = \prod_{j=1}^{n-1} ((d-1)j+1).$$
Let us first define a multivariate generating function that also incorporates the tree functional $F$ and its toll function $f$. Specifically, we set
$$Y(x,a,b) = \sum_{T} \frac{x^{|T|}}{|T|!}  e^{aF(T)-bf(T)}.$$
In view of the recursion satisfied by $F$,~\eqref{eq:diffeq2} becomes
$$\frac{\partial}{\partial x} Y(x,a,a) = \sum_{T}\frac{x^{|T|-1}}{(|T|-1)!} e^{a(F(T)-f(T))} = \Phi(Y(x,a,0)), \qquad Y(0,a,b) = 0.$$
We set
$$Z(x,a,b) = 1 + Y(xe^{-a\mu},a,b) =1 + \sum_{T} \frac{x^{|T|}}{|T|!} e^{aF(T)-a\mu |T|-bf(T)},$$
where $\mu$ will be determined later, so that
$$\frac{\partial}{\partial x} Z(x,a,a) 
= e^{-a\mu} \Phi(Y(xe^{-a\mu},a,0)) = e^{-a\mu} \Phi(Z(x,a,0)-1) = e^{-a\mu} Z(x,a,0)^d.$$
Note that
$$M_n(a) = \frac{[x^n] Z(x,a,0)}{[x^n] Z(x,0,0)} = \frac{n! [x^n] Z(x,a,0)}{Y_n}$$
is the moment generating function for the random variable $F(T_n)-\mu|T_n| = F(T_n) - \mu n$ when a random $d$-ary increasing tree $T_n$ with $n$ vertices is generated. Its derivatives with respect to $a$, evaluated at $0$, yield the moments.

Let the $r$-th derivative of $Z$ with respect to $a$ be denoted by $Z^{(r)}(x,a,b)$. Our first goal is to determine a differential equation for the function $Z^{(r)}(x,0,0)$. This is done by means of Fa\`a di Bruno's formula. First, we need some further notation regarding integer partitions: we represent partitions of a positive integer $r$ as sequences $\ell = (\ell_1,\ell_2,\ldots)$, where $\ell_j$ denotes the multiplicity of $j$. Thus $\ell$ is a partition of $r$ if $\sum_j j \ell_j = r$. The set of all partitions of $r$ is denoted by $\mathcal{P}(r)$, and we write $|\ell| = \ell_1 +\ell_2 + \cdots$ for the total number of parts in the partition $\ell$.

\begin{lemma}\label{lem:diffeq}
The function $Z^{(r)}(x,0,0)$ satisfies the differential equation
\begin{multline}\label{eq:zr_recursion}
\frac{\partial}{\partial x} \Big( Z(x,0,0)^{-d} Z^{(r)}(x,0,0) \Big) \\
= -Z(x,0,0)^{-d} H_r(x) +\sum_{s=0}^r \binom{r}{s} (-\mu)^{r-s} s! \sum_{\substack{\ell \in \mathcal{P}(s) \\ \ell_r \neq 1}}
\frac{d!}{(d-|\ell|)!} \prod_{j \geq 1} \frac{1}{\ell_j! j!^{\ell_j}} \Big( \frac{Z^{(j)}(x,0,0)}{Z(x,0,0)} \Big)^{\ell_j}\,,
\end{multline}
where
$$H_r(x) = \sum_{s=1}^r \binom{r}{s} \sum_{T} \frac{x^{|T|-1}}{(|T|-1)!} (F(T)-\mu|T|)^{r-s} (-f(T))^s.$$
\end{lemma}

Note that at this stage, $H_r(x)$ is only considered as a formal power series, convergence is not taken into account. We first analyse this differential equation in the special cases $r=1$ and $r=2$ corresponding to mean and variance before we move on to the central limit theorem.

\section{Mean and variance}\label{sec:mv}

Let us now determine mean and variance of $F(T_n)$. Since the values of the toll function $f(T)$ for $|T|>n$ will not affect the distribution  of $F(T_n)$, we can assume in this section that $f(T)=0$ for $|T|>n$. This means in particular that the functions $H_r(x)$ also depend on $n$, so we write $H_r^{(n)}(x)$ to emphasize this dependence. For $r=1$, Equation \eqref{eq:zr_recursion} becomes
$$\frac{\partial}{\partial x} \Big( Z(x,0,0)^{-d} Z^{(1)}(x,0,0) \Big) = -Z(x,0,0)^{-d} H_1^{(n)}(x)  - \mu,$$
so 
$$Z^{(1)}(x,0,0) = Z(x,0,0)^d \int_0^x \Big( -Z(w,0,0)^{-d} H_1^{(n)}(w)  - \mu \Big)\,dw,$$
where
\begin{equation}\label{eq:H1}
H_1^{(n)}(x)=-\sum_{|T|\leq n}\frac{x^{|T|-1}}{(|T|-1)!}f(T). 
\end{equation}
If we choose $\mu = \mu^{(n)}$ in such a way that
$$\mu^{(n)} = -(d-1) \int_0^{1/(d-1)}  Z(w,0,0)^{-d} H_1^{(n)}(w) \,dw,$$
then we can write

\begin{equation}\label{z1rep}
Z^{(1)} (x,0,0)=\frac{\mu^{(n)}}{d-1}Z(x,0,0)+R(x),
\end{equation}
where 
$$
R(x)=Z(x,0,0)^{d}\int_x^{1/(d-1)}  Z(w,0,0)^{-d} H_1^{(n)}(w)  \,dw.
$$
The first term on the right side of~\eqref{z1rep} contributes $\frac{\mu^{(n)}}{d-1}$ to the mean, so it suffices to determine the contribution from $R(x)$. Note that $R(x)$ is a polynomial of degree $n$ whose coefficients can be computed explicitly using the following lemma:
\begin{lemma}\label{lem2}
If 
$$
P(x)=\sum_{k=0}^{n-1}a_kx^k
 \ \text{ and } \ 
Q(x)=(1-x)^{-\beta}\int_{x}^1(1-w)^{\beta}P(w)dw,
$$
then $Q(x)$ is a polynomial of degree $n$ with 
\begin{align*}
[x^m]Q(x)
& =-\frac{a_{m-1}}{m+\beta}+\sum_{k=m}^{n-1} \binom{\beta+m-1}{m} \cdot \frac{\Gamma(\beta+1)k! a_k}{\Gamma(\beta+k+2)}\\
& =\Oh\left(\frac{|a_{m-1}|}{m}+m^{\beta-1}\sum_{k=m}^{n-1}k^{-\beta-1}|a_k |\right).
\end{align*}
\end{lemma}

This lemma gives us in particular an expression for $[x^n]R(x)$, since
\begin{align*}
[x^n]R(x) &= [x^n] (1-(d-1)x)^{-d/(d-1)} \int_x^{1/(d-1)} (1-(d-1)w)^{d/(d-1)} H_1^{(n)}(w)\,dw \\
&= (d-1)^{n-1} [x^n] (1-x)^{-d/(d-1)} \int_{x}^{1} (1-u)^{d/(d-1)} H_1^{(n)} \Big( \frac{u}{d-1} \Big)\, du.
\end{align*}
Evaluating the integral in the expression for $\mu^{(n)}$ explicitly gives us
\begin{align*}
\mu^{(n)} 
&= d(d-1) \sum_{m \leq n} \frac{1}{((d-1)m+1)((d-1)m+d) Y_m} \sum_{|T| = m} f(T) \\
&= d(d-1) \sum_{m \leq n} \frac{\E(f(T_m))}{((d-1)m+1)((d-1)m+d)}.
\end{align*}
Putting everything together, we arrive at an explicit formula for the mean:
\begin{align*}
\E(F(T_n)) &= \mu^{(n)} n + \frac{\mu^{(n)}}{d-1}+ \frac{n! [x^n] R(x)}{Y_n} \\
&= (d(d-1)n + d) \sum_{m \leq n} \frac{\E(f(T_m))}{((d-1)m+1)((d-1)m+d)} + \frac{n}{(n+d/(d-1))Y_n}  \sum_{|T| = n} f(T) \\
&= (d(d-1)n + d) \sum_{m < n} \frac{\E(f(T_m))}{((d-1)m+1)((d-1)m+d)} + \E(f(T_n)).
\end{align*}
If we complete the series and make use of conditions (C1) and (C2), we arrive exactly at the desired asymptotic formula for the mean in Theorem~\ref{thm:main}. The variance, which is obtained by using Equation \eqref{eq:zr_recursion} for $r=2$, can be treated in a similar fashion. Without going into detail, under our conditions (C1) and (C2), we can show that 
$$
\mathrm{Var}(F(T_n))=c^{(n)}n+o(n),
$$ 
where $c^{(n)}$ is a truncated double series which converges to the constant $\sigma^2$ in Theorem \ref{thm:main}  as $n\to\infty$.

\section{The central limit theorem}

We first consider the case that $f(T)$ has finite support. Conditions (C1) and (C2) are then automatically satisfied, hence the results in the previous section for the mean and variance are valid in this case as well. For the central limit theorem, we also need higher moments, for which we have the following statement.

\begin{lemma}\label{lem:finite}
If the toll function $f$ has finite support, i.e. there exists a constant $K$ such that $f(T) = 0$ whenever $|T| > K$, then the centred moments of the functional $F$ are asymptotically given by
$$\E((F(T_n)-\mu n)^r) = \begin{cases} (r-1)!! \sigma^{r} n^{r/2} + O\left(n^{r/2-1}\right) & \text{$r$ even,} \\ O\left(n^{(r-1)/2}\right) & \text{$r$ odd.} \end{cases}$$
Here, $\mu$ and $\sigma$ are as in Theorem~\ref{thm:main}. Consequently, if $\sigma \neq 0$, then the renormalised random variable 
$$\frac{F(T_n)-\mu n}{\sqrt{\sigma^2 n}}$$
converges weakly to a standard normal distribution.
\end{lemma}

The key to proving this lemma is the fact that the functions $H_r$ are now given by finite sums and therefore trivially represent entire functions. This enables us to apply singularity analysis to the functions $Z^{(r)}(x,0,0)$ for arbitrary $r$, which yields the asymptotics of the centred moments.

To deal with toll functions that are not finitely supported, we employ a trick that was already used in \cite{holmgren2015limit,janson2016asymptotic}: we approximate them by truncated versions to which we can apply Lemma~\ref{lem:finite}. This approach is based on the following simple yet general lemma.
\begin{lemma}\label{lem:0}
If $(X_n)_{n\geq 1}$ and $(W_{m,n})_{m,n\geq 1}$ are sequences of centred random variables such that  
\begin{itemize}
\item $W_{m,n}\overset{d}{\to}_n W_m$, and  $ W_m\overset{d}{\to}_m W,$ where $W$ has a continuous distribution function,
\item $\mathrm{Var}(X_n-W_{m,n})\to_n\gamma_m^2$ and $\gamma_m\to_m0$,
\end{itemize}
then 
$
X_n\overset{d}{\to}_n W.
$
\end{lemma}

We return to additive functionals and assume that the toll function $f(T)$ satisfies conditions (C1) and (C2). For every positive integer $m$, consider the truncated toll function $f_m$ and the corresponding function $F$:
$$
f_m(T)= \begin{cases} f(T) & |T| \leq m, \\ 0 & \text{otherwise,} \end{cases}\qquad\text{and}\qquad F_m(T)=\sum_{S \in \mathcal{F}(T)}f_m(S) = \sum_{S \in \mathcal{F}(T), |S|  \leq m}f(S).
$$

From Section \ref{sec:mv}, we know that the mean and variance of $F_m(T)$ have the asymptotic estimates
$$
\mathbb{E}(F_m(T))=\mu_mn+\frac{\mu_m}{d-1}+o(1) \ \ \text{ and } \ \ \mathrm{Var}(F_m(T))= \sigma^2_mn+o(n)
$$
as $n\to\infty$. Furthermore, for each $m$, if $\sigma_m^2\neq 0$ then $F_m(T)$ satisfies the central limit theorem, and  $\mu_m\to \mu$ and $\sigma^2_m\to\sigma$ as $m\to\infty$.  On the other hand, the functional $F(T)-F_m(T)$ is also additive with toll function  $f(T)-f_m(T)$. The conditions (C1) and (C2) are both satisfied by the latter toll function, so  from the asymptotic formula for the variance we know that
$$
\gamma^2_m=\lim_{n\to\infty}\frac{\mathrm{Var}(F(T_n)-F_m(T_n))}{n} \to_m 0
$$
under the conditions on the toll function $f$. Hence, Lemma \ref{lem:0} applies to the sequences
$$
W_{m,n}=\frac{F_m(T_n)-\mathbb{E}(F_m(T_n))}{n} \, \text{ and } \, X_n=\frac{F(T_n)-\mathbb{E}(F(T_n))}{n},
$$
which proves Theorem~\ref{thm:main} for arbitrary toll functions $f$ that satisfy (C1) and (C2). 

\section{Some applications}

\subsection{Fringe subtrees of given size and occurrences of specific fringe subtrees}

The simplest example of a toll function is perhaps the indicator function of a specific tree $S$:
$$f(T) = \begin{cases} 1 & T = S, \\ 0 & \text{otherwise.} \end{cases}$$
The associated additive functional is simply the number of occurrences of $S$ on the fringe of a random tree: by an occurrence of $S$, we mean a fringe subtree that is isomorphic to $S$ (including the relative order of the labels).

In this case, we obtain a central limit theorem with mean and variance only depending on the size of $S$: if $S$ has $k$ vertices, then
$$\mu = \frac{d-1}{\prod_{j=1}^{k} ((d-1)j+d)}\qquad \text{and} \qquad \sigma^2 = -\mu^2 \Big( 2k + \frac{1}{d-1} \Big) + \mu +  \frac{d(d-1)^{1-2k}}{(k-1)!^2} \int_0^1 \phi_k(x)^2\,dx.$$
A closely related functional is the number of fringe subtrees of some given size $k$ (equivalently, the number of vertices with exactly $k-1$ descendants). In particular, the special case $k=1$ corresponds to the number of leaves. Here, the toll function is given by
$$f(T) = \begin{cases} 1 & |T| = k, \\ 0 & \text{otherwise,} \end{cases}$$
and we obtain a central limit theorem with
$$\mu = \frac{d(d-1)}{((d-1)k+d)((d-1)k+1)} \quad \text{and} \quad \sigma^2 = -\mu^2 \Big( 2k + \frac{1}{d-1} \Big) + \mu +  \frac{d(d-1)^{1-2k}Y_k^2}{(k-1)!^2} \int_0^1 \phi_k(x)^2\,dx.$$
This was already shown by Fuchs \cite{fuchs2012limit}, who also considered the case that $k$ is not fixed but rather tends to infinity with the size of the tree as well.

\subsection{The number of subtrees}

The number of subtrees is already a somewhat more complicated example: for Galton-Watson trees, binary increasing trees and recursive trees, it was already studied in \cite{wagner2015central}. Here, we count all subtrees, i.e. all induced subgraphs that are again trees, not just those on the fringe. It is useful to study an auxiliary quantity first, namely the number of subtrees containing the root: we write $s(T)$ for this number. It is not difficult to see that
$$s(T) = \prod_{j=1}^k (1+s(B_j)),$$
since each subtree induces either the empty set or a subtree containing the root in each of the branches. Taking the logarithm gives us
$$\log(1+s(T)) = \sum_{j=1}^k \log(1+s(B_j))  + \log(1+s(T)^{-1}),$$
so $\log(1+s(T))$ is an additive functional with toll function $f(T) = \log(1+s(T)^{-1})$. Simple a priori estimates show that the technical conditions of our general central limit theorem are satisfied: this is because $s(T) \geq |T|$ (since every path from the root to a vertex is also a subtree), which implies that $f(T) = \Oh(|T|^{-1})$ for all $T$ (even deterministically, not just on average). Thus our main result applies to the functional $s(T)$. As it was shown in \cite{wagner2015central}, the difference between $F(T) = \log(1+s(T))$ and the logarithm of the total number of subtrees (not necessarily containing the root) is $\Oh(\log |T|)$, so the central limit theorem remains correct for the total number of subtrees.

\subsection{The size of the automorphism group}

An important motivating example for this paper is the size of the automorphism group. In their article \cite{bona2009isomorphism}, B\'ona and Flajolet proved, motivated by questions in phylogenetics, that the logarithm of the size of the automorphism group of uniformly random binary trees is asymptotically normally distributed (they proved this limit law for the number of nodes for which the two branches are isomorphic, which is equivalent). Here, we obtain an analogous statement for $d$-ary increasing trees. We remark that binary increasing trees are also essentially equivalent to the Yule-Harding model (as opposed to the uniform model) of phylogenetics \cite[Section 2.5]{semple2003phylogenetics}.

As it was mentioned in the introduction, the relevant toll function is $f(T) = \log (R(T))$, where $R(T)$ is the size of the symmetry group of the collection of root branches. This simplifies considerably in the case of binary trees, where we only have two branches $B_1$ and $B_2$. In this case, it follows that
$$f(T) = \begin{cases} \log 2 & \text{if } B_1 \text{ and } B_2 \text{ are isomorphic,} \\ 0 & \text{otherwise.}\end{cases}$$
As one would expect, it is very unlikely for large trees that the two branches are actually isomorphic, which is why the technical condition on the toll function is satisfied. In fact, one can show that $\E(|f(T_n)|)$ decays exponentially for binary increasing trees. We find that the number of automorphisms of a random binary increasing tree asymptotically follows a log-normal law, which parallels the aforementioned result of B\'ona and Flajolet.

The same holds more generally for $d$-ary trees, although the expected value of the toll function does not decay as quickly: in this case, the probability that two branches are isomorphic only decreases at a rate of $\Oh(|T|^{-2/(d-1)})$, which however is still sufficient.

\subsection{The number of orbits}

Two vertices of a rooted tree (or generally any graph) are said to belong to the same orbit if there exists an automorphism that maps one of the vertices to the other. The vertex set can thus be partitioned in a natural way into orbits, and the number of orbits can also be regarded as an additive functional. Let us illustrate this for binary increasing trees: if the two branches $B_1$ and $B_2$ of $T$ are isomorphic, then $F(T)=F(B_1)+1=F(B_2)+1$, otherwise, $F(T)=F(B_1)+F(B_2)+1$. Hence, the toll function is given by
$$f(T) = \begin{cases} 1-F(B_1) & \text{if } B_1 \text{ and } B_2 \text{ are isomorphic,} \\ 1 & \text{otherwise.}\end{cases}$$
Our technical conditions (C1) and (C2) are not completely satisfied in this examples, but it is possible to work around that. First, the expected value of $|f(T_n)|$ does not tend to $0$, but the expected value of $|f(T_n)-1|$ does, cf. Remark~\ref{rem:constant}. Second, $f(T)$ is not bounded, but this condition can be replaced by the fact that $f(T)$ and $F(T)$ are both $O(|T|)$. Again, everything also remains valid for $d$-ary increasing trees, although the details are somewhat more intricate.

\bibliographystyle{abbrv}
\bibliography{increasing}
\end{document}